\documentclass[11pt,reqno]{amsart}
\usepackage{amssymb,amsmath}\newcommand{\be}{\begin{eqnarray}}
\newcommand{\ee}{\end{eqnarray}}
\newcommand{\Id}{\mbox{\rm Id}}

\newcommand{\vf}{\varphi}

\newcommand{\half}{\frac12}

\newcommand{\R}{{\mathbb R}}\newcommand{\pd}{\partial}

\newcommand{\C}{{\mathbb C}}

\newcommand{\B}{{\mathcal B}}

\newcommand{\K}{{\mathcal K}}\newcommand{\M}{{\mathcal M}}

\newcommand{\rank}{{\rm rank}}

\newcommand{\F}{\mathcal{F}}

\newtheorem{theorem}{Theorem}\newtheorem{lemma}[theorem]{Lemma}

\theoremstyle{definition}
\theoremstyle{remark}\numberwithin{equation}{section}\input epsf.sty
\begin{document}\thispagestyle{empty}

\title[Ahlfors-Beurling operator on radial functions]{{Ahlfors-Beurling operator on radial functions}}
\author{Alexander Volberg}\address{Alexander Volberg, Department of  Mathematics, Michigan State University and the University of Edinburgh.
{\tt volberg@math.msu.edu}\,\,and\,\,{\tt a.volberg@ed.ac.uk}}
\date{ 5 February 2009}
\begin{abstract}
We give a simple proof of a sharp bound of Ahlfors--Beurling operator on complex-valued radial functions. In the language of the Calculus of Variation we prove a certain estimates for stretch and twist functions. Our estimate go slightly beyond this case. This article was written in February 2009 and then delivered at a seminar talk at UW, Madison on February, 2009.
\end{abstract}
\maketitle

\section{Introduction}
\label{1}

``Everything has been thought of before, the task is to think about it again" said Goethe. We want to take another look at Ahlfors-Beurling operator $T$, it is the operator that sends
$\bar\pd f$ to $\pd f$ for smooth functions $f$ with compact support on the plane $\C$.
Here 
$$
\pd f =\frac{\pd f}{\pd z} = \bigg( \frac{\pd f}{\pd x} - i \frac{\pd f}{\pd y}\bigg)\,,\,\,\bar\pd f =\frac{\pd f}{\pd z} = \bigg( \frac{\pd f}{\pd x} +i \frac{\pd f}{\pd y}\bigg)\,.
$$
We intentionally omitted $\half$, this will not bring complications.

 This operator was much studied in the last 30 years. There are several reasons for that.

a) Operator $T$ and its multidimensional analogs play an important part in the theory of quasiregular mappings.

b) Attempts to calculate (estimate) the norm of $T$ are closely related to important conjectures in the Calculus of Variation: Morrey's conjecture of 1952 and Sverak's conjecture of 1992. Morrey's conjecture states that ``rank one convex functions are not necessarily quasiconvex", so in essence it asks for a series of counterexamples, of rank one convex functions that are not quasiconvex. Sverak's conjecture asks about a {\it concrete} rank one convex function whether it is quasiconvex.

c) There is a deep connection of Ahlfors-Beurling operator to stochastic calculus and stochastic optimal control.

Saying all that let us state several very innocent looking problems.

\section{Some problems}
\label{problems}

We mostly follow in this section the exposition of A. Baernstein--S. Montgomery-Smith \cite{BaMS}.
 
 Define a function $L:\C^2\rightarrow \R$ as follows
 $$
 L(z,w)=\begin{cases} |z|^2-|w|^2\,,\,\,\, \text{if}\,\,\,  |z|+|w|\le 1 \,,\\
 2|z|-1\,,\,\,\,\text{if}\,\,\, |z|+|w|>1\,.\end{cases}
 $$

\bigskip

{\it Sverak's problem}: Let $f \in C_0^{\infty}(\C)$. Is it true  that
\begin{equation}
\label{Svconj}
\int_{\C} L(\bar\pd f, \pd f) dxdy \ge 0\,?
\end{equation}

We can restate this problem in the language of quasiconvex functions of matrix argument. Then we will explain Morrey's problem. 

Let $M(m,n)$ be the set of all $m\times n$ matrices with real entries. A function $\Psi: M(m,n)\rightarrow \R$ is called {\it rank one convex} if $t\rightarrow \Psi (A +tB)$ is convex function for any $B\in M(m,n)$ that has rank $1$ and any $A\in M(m,n)$.

Function $\Psi$ is called {\it quasiconvex} if it is locally integrable, and for each $A\in M(m,n)$ and each bounded domain $\Omega\subset \R^n$ and each smooth function $f:\R^n\rightarrow \R^m$ one has
\begin{equation}
\label{qconv}
\frac{1}{|\Omega|}\int_{\Omega} \Psi (A + Df) dx \ge \Psi(A)\,.
\end{equation}
Here $Df$ is the Jacobi matrix of the map $f$.

For $n=1$ or $ m=1$ quasiconvexity is equivalent to convexity (which of course is equivalent for this case to rank one convexity). Always convexity implies quasiconvexity that implies rank one convexity.

\bigskip

{\it Morrey's problem}: If $m>1, n>1$ rank one convexity does not imply quasiconvexity. This was conjectured by Morrey in 1952 in \cite{Mo}. Sverak \cite{Sv2} proved that problem if $m>2$. If $m=2$ this is still open even in the case $n=2$. Morrey's problem enjoyed a lot of attention in the last 57 years.  

\bigskip

We can translate easily Sverak's problem to this language (this is how it appeared in the first place).

In fact,  for $A\in M(2,2)$, $A=\begin{bmatrix} a\,& b\\c\,&d\end{bmatrix}$. Put $z= a-d +i(b+c)$,
$w=a+d +i (c-b)$. We see that
\begin{equation}
\label{detdef}
\Psi(A) := L(z,w) =\begin{cases} -4\det (A)\,,\,\,\,\text{if}\,\,\, (|A|_2^2 - 2\det A)^{\frac12} + (|A|_2^2 - 2\det A)^{\frac12} \le 1\,,\\ 2(|A|_2^2 - 2\det A)^{\frac12} -1\,,\,\,\,\text{otherwise}\,.\end{cases}
\end{equation}

This function is rank one convex on $M(2,2)$.  A very simple proof is borrowed from \cite{BaMS}.
We fix $A, B\in M(2,2)\,, \rank (B)=1$. Let $(z, w)$ corresponds to $A$ and $(Z, W)$ to $B$. The fact that $\rank (B)=1$ means that the map $\zeta\rightarrow Z\zeta +W \bar\zeta$ maps the plane to the line, so $|Z|=|W|$. Then $|z+tZ|^2 -|w+tW|^2 = a + tb$ for some $a,b\in \R$---there is no quadratic term. Also
$\Psi(A+tB) = |z+tZ|^2 -|w+tW|^2 = a + tb$ if and only if $|z+tZ| +|w+tW| \le 1$. As all $z,Z, w, W$ is fixed and $t\rightarrow |\alpha +t \beta|$ is convex for any complex $\alpha, \beta$, we conclude that
$\{t\in\R: |z+tZ| +|w+tW| \le 1\}$ is an interval (may be empty). On the other hand outside of this interval $\Psi(A+tB) = 2|z+tZ|-1$, that is a convex function. Now continuity of $\Psi(A+tB)$ implies that it is convex.

\bigskip

Let $f:\R^2\rightarrow\R^2$, smooth, and with compact support. Write $f=u+iv$, $Df =  \begin{bmatrix} u_x \,& u_y\\v_x\,&v_y\end{bmatrix}$. Then
using the notations above $z= u_x - v_y +i( v_x +u_y) =\bar\pd f$, $w= u_x + v_y +i( v_x-u_y)=\pd f$.

This means that $\Psi (Df) = L(\bar\pd f, \pd f)$ and Sverak's conjecture states that
$$
\int \Psi(Df) dxdy \ge \Psi(0) =0\,.
$$ 

In other words, \eqref{qconv} means that $\Psi$ from \eqref{detdef} is quasiconvex at $A=0$.  We conclude that Sverak's conjecture is not true, then $\Psi$ gives an example of rank one convex function which is not quasiconvex. This would solve Morrey's conjecture which exactly asks for such an example for the case $n=2, m=2$. 

\bigskip

However, \eqref{qconv} is probably true. Everybody who worked with these questions believes in it. We will explain this belief.

\section{Consequences of Sverak's inequality \eqref{qconv}.}
\label{cons}

In what follows 
$$
p^* = \max (p, \frac{p}{p-1}) =\max (p, p')\,.
$$

Here is one other function on $M(2,2)$ which is rank one convex but for which it is unknown whether it is quasiconvex. It is also on $M(2,2)$. Several such functions are discussed in \cite{Sv1}, \cite{Sv2}, but the function $\Psi$ above and $\Psi_p$ below are especially important for us.

$$
\Psi_p(A) = 
$$
$$
((p^* -1) (|A|_2^2 -2\det A)^{\frac12} - (|A|_2^2 +2\det A)^{\frac12}) ((|A|_2^2 -2\det A)^{\frac12}+(|A|_2^2 +2\det A)^{\frac12})^{p-1}\,.
$$

Repeat our correspondence between real matrices $M(2,2)$ and $\C^2$: for $A\in M(2,2)$, $A=\begin{bmatrix} a\,& b\\c\,&d\end{bmatrix}$, put $z= a-d +i(b+c)$,
$w=a+d +i (c-b)$. We see that
\begin{equation}
\label{detdefp}
\Psi_p(A) := L_p(z,w) = ((p^*-1) |z|-|w|) (|z|+|w|)^{p-1}\,.
\end{equation}

See now e. g. \cite{BaMS} for
$$
L_p(z,w) = \frac2{p(2-p)}\int_0^{\infty} t^{p-1} L(\frac{z}t, \frac{w}t) dt\,,\,\,\,\text{if}\,\, 1<p<2\,.
$$
Obviously, for any $z, w, Z, W, |Z|=|W|$ the function $t\rightarrow L_p(z+tZ, w+tW)$ is convex because of the formula and because we just proved such a property for  $L$.

Then, automatically, 
$$
\Psi_p(A) = \int_0^{\infty} t^{p-1} \Psi(\frac{A}t) dt\,,\,\,\,\text{if}\,\, 1<p<2\,.
$$

And then $\Psi_p$ is a rank one convex function in an obvious way, if $1<p<2$.

But for $2< p<\infty$ another formula holds (see again \cite{BaMS}): put
$$
M(z,w) = L(z,w)- (|z|^2-|w|^2) = (|w|^2 -(|z|-1)^2){\bf 1}_{|z|+|w|>1}\,.
$$
Obviously, for any $z, w, Z, W, |Z|=|W|$ the function $t\rightarrow M(z+tZ, w+tW)$ is convex because we subtract the linear term $a+bt$ from $L(z+tZ, w+tW)$.

Then 
$$
\M_p(z,w) = \frac2{p(p-1)(p-2)}\int_0^{\infty} t^{p-1} M(\frac{z}t, \frac{w}t) dt\,,\,\,\,\text{if}\,\, 2<p<\infty\,,
$$
is such that for any $z, w, Z, W, |Z|=|W|$ the function $t\rightarrow M_p(z+tZ, w+tW)$ is convex.

And, automatically, 
$$
\Psi_p(A) = \int_0^{\infty} t^{p-1} (\Psi(\frac{A}t) + \frac4{t^2}\det A)dt\,,\,\,\,\text{if}\,\, 2<p<\infty\,,
$$
is a rank one convex function on $M(2,2)$.

\bigskip

{\it Banuelos-Wang problem}:
Is it true that for any smooth function with compact support on $\C$
\begin{equation}
\label{BWpr}
\int_{\C} L_p(\bar\pd f, \pd f) dx dy \ge 0\,?
\end{equation}

If \eqref{BWpr} were {\it not} true we would have that $\Psi_p$ is {\it not} quasiconvex at $A=0$ and Morrey's problem would be solved in the remaining case.

If \eqref{BWpr} were true than we would have solved Iwaniec's problem of 1982.

{\it Iwaniec's problem}: Ahlfors-Beurling operator T which sends $\bar\pd f$ to $\pd f$ has norm $p^*-1$.
 Essentially it is the following inequality for all $f\in C_0^{\infty}(\C)$:
 \begin{equation}
\label{Iwpr}
\int_{C} |\pd f|^p dx dy \le (p^*-1)^p\int_{C} |\bar\pd f|^p dx dy  \,?
\end{equation}
In equivalent form \eqref{Iwpr} is stated as follows

\begin{equation}
\label{IwprT}
\int_{C} |Tf|^p dx dy \le (p^*-1)^p\int_{C} |f|^p dx dy\,,\,\,\,\text{for all}\,\,\, f\in C_0^{\infty}(\C)  \,?
\end{equation}

In fact, \eqref{BWpr} $\Rightarrow$ \eqref{Iwpr} follows from a pioneering research of Burkholder, who in \cite{Bu1}, \cite{Bu3}, p. 77,
noticed that

\begin{equation}
\label{Bur1}
p\bigg(1-\frac1{p^*}\bigg)^{p-1} L_p(z,w) \leq (p^*-1)^p |z|^p - |w|^p\,.
\end{equation}
Now it is clear why  \eqref{BWpr} implies \eqref{Iwpr}. 

\bigskip

{\bf Remark.} What is subtle and interesting is the whole theory of inequalities of the type like Burkholder's inequality \eqref{Bur1}. This is actually the essence of the so-called Bellman function approach. The literature is now extensive, and it relates \eqref{Bur1} to Monge-Amp\`ere equation and stochastic control, see e. g. Slavin-Stokolos' paper \cite{SlSt} or Vasyunin and Volberg \cite{VaVo2}.

\bigskip

Sverak's conjecture \eqref{Svconj} and, as a result, Banuelos-Wang's conjecture \eqref{BWpr} were proved in the paper of Baernstein and Montgomery-Smith \cite{BaMS} in the case of so-called ``stretch functions" $f$. A stretch function (in our notations, which differ slightly from those in \cite{BaMS}) is a function of the form
$$
f(re^{i\theta}) = g(r) e^{-i\theta}\,,
$$
where $g$ is a smooth function on $\R_+$, $g(0)=g(\infty)=0$, and $g\ge 0$. We will call such $g$'s {\bf stretches}.

A straightforward calculation  shows:
\begin{equation}
\label{d}
\bar\pd f = g'(r) +\frac{g(r)}{r}\,,\,\,\, \pd f = e^{-2i\theta} \bigg(g'(r) - \frac{g(r)}{r}\bigg)\,.
\end{equation}

So in \cite{BaMS} it is proved that for any stretch $g$ (in particular, $g$ must be non-negative)
\begin{equation}
\label{g}
\int_0^{\infty}|g'(r) -\frac{g(r)}{r}|^p \,rdr\le (p^*-1)^p\int_0^{\infty} |g'(r) + \frac{g(r)}{r}|^p \,rdr\,, 1<p<\infty\,.
\end{equation}

Let us change the variable: 
\begin{equation}
\label{beta}
\beta(\rho) := \half\bigg(g'(\sqrt{\rho}) + \frac{g(\sqrt\rho)}{\sqrt\rho}\bigg)\,,\rho\in \R_+\,.
\end{equation}
If we introduce Hardy operator:
$$
H \beta(u) = \frac1u\int_0^u \beta(s)\,ds
$$ on locally integrable functions on $[0, \infty)$, we an invert \eqref{beta} for any $g\in C_0^{\infty}(\R_+)$:

\begin{equation}
\label{invbeta}
g(\sqrt\rho) = \sqrt\rho (H\beta)(\rho)\,.
\end{equation}
In fact, if we define $g_1(\sqrt\rho)=\sqrt\rho (H\beta)(\rho)$, where $\beta$ is from \eqref{beta}, we get
that both $g, g_1$ satisfy \eqref{beta} (an easy calculation for $g_1$). Let $g_2= g-g_1$. Then
$$
g_2'(\sqrt{\rho}) + \frac{g_2(\sqrt\rho)}{\sqrt\rho} =0\,,\,\,\forall \rho\in\R_+\,.
$$
Consider $f(re^{i\theta}) := g_2(r) e^{-i\theta}$. The previous formula and our previous calculation of $\bar\pd f=
g_2'(r) +\frac{g_2(r)}{r}$ shows that $\bar\pd f=0$. Function $f$ is entire and vanishes at infinity, this $|g_2|=  |f|=0$. 

On the other hand, \eqref{invbeta} implies
$$
\frac12   \bigg(g'(\sqrt{\rho}) - \frac{g(\sqrt\rho)}{\sqrt\rho}\bigg) =\rho (H\beta)'(\rho) =\beta(\rho)-(H\beta)(\rho)\,.
$$

 So \cite{BaMS} proves that for all $\beta=g'(\sqrt{\rho}) + \frac{g(\sqrt\rho)}{\sqrt\rho}$, where $g$ is a {\bf stretch} one has
 \begin{equation}
 \label{H_I}
 \int_0^{\infty} |(H-I) \beta(\rho)|^p\,d\rho \leq (p^*-1)^p \int_0^{\infty} |\beta(\rho)|^p \, d\rho\,.
 \end{equation}  
 
 In particular, \eqref{H_I} holds for all $\beta$ such that $H\beta\ge 0, H\beta\in C_0^{\infty}(\R_+)$.
In the paper of Banuelos and Janakiraman \cite{BaJa2} it was observed that such $\beta$'s are dense in $L^p_{\text{real}}(\R_+)$.
Therefore, \eqref{H_I} means

\begin{equation}
 \label{normH_Ireal}
 \|H-I\|_{L^p_{\text{real}}(\R_+)\rightarrow L^p_{\text{real}}(\R_+)} \le p^*-1\,.
 \end{equation}
 
 It is interesting to compare this with classical Hardy's inequality:
 
 \begin{equation}
 \label{normH}
 \|H\|_{L^p_{\text{real}}(\R_+)\rightarrow L^p_{\text{real}}(\R_+)} \le p^*\,, 1<p\le 2\,.
 \end{equation}
 
 And both results are sharp for $1<p\le 2$: 
 $$
 \|H-I\|_{L^p_{\text{real}}(\R_+)\rightarrow L^p_{\text{real}}(\R_+)} =  p^*-1\,,\,\,\|H\|_{L^p}= p^*\,,\,1<p\le 2\,.
 $$
 
 But the word ``real" can  betaken out by Marcinkiewicz--Zygmund's lemma (which says that the operator with the real kernel will have the same norm on complex-valued and real-valued functions), and we come to
 \begin{equation}
 \label{normH_I}
 \|H-I\|_{L^p(\R_+)\rightarrow L^p(\R_+)} \le p^*-1\,.
 \end{equation}
 
 From this we can easily conclude that
 for {\bf complex stretches} $g$ (smooth, compactly supported) inequality \eqref{g} holds.
 
 Therefore,
 \begin{equation}
\label{IwprT1}
\int_{C} |Tf|^p dx dy \le (p^*-1)^p\int_{C} |f|^p dx dy\,,\,\,\,\text{for all}\,\,\, f\in C_0^{\infty}(\C)\,,f(z)=f(|z|)\,, f:\C\rightarrow\C  
\end{equation}
holds for all complex valued radial $f$. This was a question in \cite{BaJa2}. 

\bigskip

We can also show how to do the estimate on complex-valued radial  functions using Bellman function techniques. The interest of that is in the fact that we can go a bit beyond the radial functions. Otherwise Marcinkiewicz--Zygmund lemma is enough. But the advantage of the method  below is that it is applicable to other situations. It also illustrate how genuinely convex functions can sometimes be involved in a rather sophisticated way in proving {\it quasiconvexity} statements.

\section{Bellman function and Ahlfors-Beurling operator on radial functions.}
\label{Bfsection}

The kernel of $T$ is $K(z) = \frac1{\pi}\frac1{z^2}= :e^{-2i\theta}k(r)\,,z=re^{i\theta}$. So for radial $g$
$$
Tg(\rho e^{i\vf}) = \int_{\C} k(z)g(|z-\rho e^{i\vf}|) dA(z) =e^{-2i\vf}\int\int e^{-2i\psi}k(r) g(|re^{i\psi} -\rho|)\,rdrd\psi =
$$
$$
e^{-2i\vf}\int K(w+\rho) g(|w|) dA(w)= e^{-2i\vf}\int (\int_0^{2\pi} K(|w|e^{it} +\rho)\,dt) g(|w|) dA(w)\,.
$$

\subsection{Symmetrization}
\label{symm}

If we denote $n(\rho,r)=\int_0^{2\pi} K(re^{it} +\rho)\,dt $, and $Ng(\rho):= \int_0^{\infty} n(\rho,r) g(r)\,rdr$ we get
$$
Tg(\rho e^{i\vf}) = e^{-2i\vf}\int_0^{\infty} n(\rho,r) g(r)\,rdr=: e^{-2i\vf}Ng(\rho)\,.
$$

Hence to check the norm of  
$Tg(\rho e^{i\vf})$ in $L^p$ we can take a function $f\in L^{p'}(\C)$, write {\it bilinear} form
$$
(f, Tg)=\int_{\C} f Tg \, dA = \int (\int e^{-2i\psi} f(re^{i\psi}) d\psi)Ng(r)\,rdr\,.
$$

Let us notice that the family $\F$ of functions having the form
$$
f(re^{i\theta}) = \sum_{k=-N}^N e^{-ik\theta} f_k(r)\,,
$$
where $f_k$ are smooth compactly supported functions, give us a dense family in $L^{p'}(\C),\,, 1<p<\infty$. Also let us call $e^{-ik\theta} f_k(r)$ a $k$-mode function.
The set of $k$-modes is called $\F_k$. 
Continuing the last formula we write
$$
(f, Tg)=\int_{\C} f Tg \, dA =
$$
\begin{equation}
\label{samemode}
 2\pi\int f_{-2}(r) Ng(r)\,rdr= \int_{\C} f_{-2}(|z|) Ng(|z|) \, dA(z)= (e^{2i\theta}f_{-2}(|z|), Tg)\,.
\end{equation}
Let us notice that projection $\Pi_k:\F\rightarrow \F_k$ has norm at most $1$ in any $L^p$.
In fact,
Let $R_{\vf}$ is a rotation of $\C$ by $\vf$. Then
$$
f_k(|z|) = \frac1{2\pi} \int_0^{2\pi} e^{ik\vf}f(R_{\vf}z)\,d\vf\,.
$$
So projection $\Pi_k$ is just the averaging-type operator, and thus has norm at most $1$.

Conclusion:
to estimate $\|Tg\|_p$, $g\in \F_0$, it is sufficient to estimate the bilinear form $|(f,Tg)|$ only for $f\in \F_{-2}$ (and in the unit ball of $L^{p'}(\C)$).
We proved actually the following
\begin{lemma}
\label{k_2}
For $g\in \F_k$,
$$
\|Tg\|_p =\sup_{f\in L^{p'}\cap \F_{k-2}\,,\|f\|_{p'}\le 1} |(f, Tg)|\,.
$$
\end{lemma}

We actually repeated also the following well-known  simple calculation.

\begin{lemma}
\label{symmkernel}
Let a complex valued kernel $K(re^{i\theta}) = e^{-il\theta} k(r)$. Let $\K f:= K\star f$ be a convolution operator.
Then it maps $\F_k$ to $F_{k-l}$ and for every $g= e^{-ik\theta}g_k(r)\in F_k$ we have
$$
\K g(\rho^{i\vf}) = e^{-i(k-l)\vf}\int_0^{\infty} N_{k}(\rho,r) g_k(r)\,rdr\,,
$$
where
$$
N_k(\rho,r) := \int_0^{2\pi} K(re^{it} +\rho) e^{-ikt}\, dt\,.
$$
\end{lemma}

\bigskip

For $\K=T$ one can compute the kernel of $N_m$:
$$
\half N_m(t,x) = x\delta_x - (m+1) \frac{1}{t^{m+2}} x^m{\bf 1}_{[0,t]}(x)\,.
$$
It is  not very nice, but let us denote by $h:\R_+\rightarrow \R_+$ the map $h(t)=t^2$. Then
the operator $\Lambda_m$ (see \cite{BaJa2} for this) is 
$$
\Lambda_m g(u) = g(u) - (m+1) \frac1{u^{\frac{m+2}2}}\int_0^u  v^{\frac{m}2}\,g(v)\,dv\,.
$$

For $m=0$ this is $\Lambda_0= \Id - H$, where $H $ is Hardy's averaging operator on half-axis:
$$
Hg(u) := \frac1u\int_0^u g(v) \,dv\,.
$$
Famous Hardy's inequality is practically equivalent to computing
$$
\|H\|_{l^p(\R_+)\rightarrow L^p(\R_+)} = p^*\,,\,\, \text{if}\,\, 1<p\le 2\,.
$$

Curiously, we can see now that the question about complex valued radial functions from \cite{BaJa2} is equivalent to
\begin{equation}
\label{HI}
\|H-\Id\|_{l^p(\R_+)\rightarrow L^p(\R_+)} \le p^*-1\,,\,\, \text{if}\,\, 1<p\le 2\,.
\end{equation}

\subsection{A Bellman function}
\label{Bf}

We will use a certain interesting convex functions on $\R^6$ and $\R^4$ to approach our ``quasiconvexity" inequality \eqref{Iwpr} for complex valued radial functions. 

Suppose we have function $\B(u,v, \xi, \eta, H,Z)$ of $6$ real variables defined in
$$
\Omega= \{|(u,v)|^p \le H\,, |(\xi,\eta)|^{p'} < Z\}\,,
$$
and satisfying

\noindent I) For an arbitrary $a\in\Omega\alpha\in \R^6$ we want to have 
$$
\langle -\frac{d^2 \B}{d a^2}\alpha,\alpha\rangle\ge 2 (\alpha_1^2+ \alpha_2^2)^{1/2}(\alpha_3^2+ \alpha_4^2)^{1/2}\,.
$$
and

\noindent II) For an arbitrary $a\in\Omega$
$$
\B(a) \le (p^*-1) \bigg(\frac{H}{p}  + \frac{Z}{p'}\bigg)\,,\,\,\text{where}\,\, p^* =\max (p,p')\,.
$$

For the sake of future convenience we prefer to work with the following transformation of $\B$ ($a=(u,v,\xi,\eta, H,Z)$):
$$
B(u,v, \xi,\eta) :=\sup_{a\in\Omega} \{B(a) - (p^*-1)\bigg(\frac{H}{p}  + \frac{Z}{p'}\bigg)\}\,.
$$

Then it is not difficult to check that this $B$ is still concave (in spite of being {\i supremum } of concave functions):

\begin{equation}
\label{d2}
-d^2B \ge 2 |(du,dv)| |(d\xi,d\eta)|\,.
\end{equation}
\begin{equation}
\label{bds}
-(p^*-1)\bigg(\frac{|(u,v)|^p}{p}  + \frac{|(\xi,\eta)|}{p'}\bigg) \le B(u,v,\xi,\eta) \le 0\,.
\end{equation}

The existence of such $\B$ was proved in \cite{PV},\cite{DV1}.

\subsection{Heat extension}
\label{heat}

Let $f,g$ be two test functions on the plane. By the same letters we denote their heat extensions into $\R_+^3$.
This is a simple lemma observed in \cite{PV}:

\begin{lemma}
\label{hext}
$$
\int_{\C} f Tg \,dA = -2 \int_{\R_+^3} (\pd_x +i\pd_y) f \cdot (\pd_x +i\pd_y)g \,dxdydt\,.
$$
\end{lemma}

Let us use  below the following notations:
$$
f=u+iv,  z_1= u_x+iu_y, z_2 = v_x+iv_y,
$$
$$
g=\xi+i\eta, \zeta_1 = \xi_x +i\xi_y, \zeta_2 = \eta_x +i\eta_y\,.
$$
 Now we can read Lemma \ref{hext} as follows:
 
 \begin{equation}
 \label{zext}
 \int fTg = -2\int_{\R_+^3} (z_1+iz_2)(\zeta_1 + i \zeta_2)\,,  |\int fTg| \le 2\int_{\R_+^3} |z_1+iz_2||\zeta_1 + i \zeta_2|\,.
 \end{equation}
 And from here we see
 \begin{equation}
 \label{zext1}
 |\int fTg| \le  2\int\bigg[\frac{|z_1+iz_2|^2 + |z_1-iz_2|^2}{2}\bigg]^{1/2}\bigg[\frac{|\zeta_1+i\zeta_2|^2 + |\zeta_1-i\zeta_2|^2}{2}\bigg]^{1/2}\,.
 \end{equation}
 
 Property \eqref{d2} of $B$ can be rewritten
 \begin{lemma}
 \label{d2z}
$$ -\langle d^2 B (z_1, z_2, \zeta_1, \zeta_2)^{T}, (z_1, z_2, \zeta_1, \zeta_2)^{T}\rangle \ge 2 [|z_1|^2 + |z_2|^2]^{1/2} [|\zeta_1|^2 + |\zeta_2|^2]^{1/2}\,.$$
 \end{lemma}
 
 This lemma gives now
 
 \begin{align*}
 \label{mix}
 -2\langle d^2 B (\frac{z_1+i z_2}{2}, \frac{z_1-i z_2}{2}, \frac{\zeta_1+i\zeta_2}{2}, \frac{\zeta_1-i\zeta_2}{2})^{T}, (\text{the same vector})^{T}\rangle \ge \\2\bigg[\frac{|z_1+iz_2|^2 + |z_1-iz_2|^2}{2}\bigg]^{1/2}\bigg[\frac{|\zeta_1+i\zeta_2|^2 + |\zeta_1-i\zeta_2|^2}{2}\bigg]^{1/2}\,.
 \end{align*}
 
 After integration and using \eqref{zext1}
 we get
 \begin{equation}
 \label{LHS}
 |\int_{\C} fTg|\le \int_{\R_+^3} LHS\,.
 \end{equation}
 
 The rest is the estimate of $\int_{\R_+^3} LHS$ from above. First of all simple algebra ($a:=(u,v,\xi,\eta)$:
 \begin{align*}
 \int_{\R_+^3} LHS = -\half\int_{\R_+^3}  \langle d^2B(a) (z_1, z_2, \zeta_1, \zeta_2)^T, (\text{the same})^T\rangle -\\ \half\int_{\R_+^3} \langle d^2B(a) (z_2, -z_1, \zeta_2, -\zeta_1)^T, (\text{the same})^T\rangle + \\ \int_{\R_+^3} \text{auxiliary terms} =: I + II + III\,.
\end{align*}

It has been proved in \cite{PV}, \cite{DV1} that (the convention is that $u,v, \xi, \eta$ are heat extensions of homonym functions on the plane)

$$
I = \half\int_{\R_+^3}\bigg(\frac{\pd}{\pd t} -\Delta\bigg) B(u,v, \xi, \eta)\,,
$$
$$
II = \half\int_{\R_+^3}\bigg(\frac{\pd}{\pd t} -\Delta\bigg) B(v,-u, \eta, -xi)\,.
$$

\bigskip

{\bf An estimate of I from above.}
Let $H$ denote the heat extension of function $|f|^{p'} = (u^2 +v^2)^{p'/2}$, $Z$ denote the heat extension of function $|g|^p = (\xi^2 +\eta^2)^{p/2}$. In $\R_+^3$ consider $\Psi(x,y,t)= B(u,v, \xi,\eta) + (p^* -1) \bigg(\frac{H}{p'}+\frac{Z}{p}\bigg)$.  Then 
$$
2I = \int_{\R_+^3}\bigg(\frac{\pd}{\pd t} -\Delta\bigg) \Psi\,,
$$
Then obviously (integration by parts)
$$
2I= \int_{\R_+^3}\frac{\pd}{\pd t}\Psi = \lim_{t\rightarrow\infty} \int_{\R^2} \Psi(\cdot,t) - \int_{\R^2} \Psi(\cdot,0)  \le \lim_{t\rightarrow\infty} \int_{\R^2} \Psi(\cdot,t)\,.
$$
Using \eqref{bds} ($B\le 0$) we get
$$
I \le  \half(p^*-1) \lim_{t\rightarrow\infty} \int_{\R^2}\Bigg(\frac{H(\cdot,t}{p'} + \frac{Z(\cdot, t}{p}\bigg)=
\half(p^*-1) \bigg(\frac{\|f\|_{p'}^{p'}}{p'} + \frac{\|g\|_p^p}{p}\bigg)\,.
$$
Similarly,
$$
II \le  
\half(p^*-1) \bigg(\frac{\|-if\|_{p'}^{p'}}{p'} + \frac{\|-ig\|_p^p}{p}\bigg)\,.
$$
So 
$$
I + II \le (p^*-1) \bigg(\frac{\|f\|_{p'}^{p'}}{p'} + \frac{\|g\|_p^p}{p}\bigg)\,.
$$

We are going to prove next that
$$
III\le 0\,.
$$
Combining we get
$|\int_{\C} fTg| \le (p^*-1) \bigg(\frac{\|f\|_{p'}^{p'}}{p'} + \frac{\|g\|_p^p}{p}\bigg)$ and the usual polarization argument proves out final statement:

$$
|\int_{\C} fTg| \le (p^*-1)\|f\|_{p'}\|g\|_p\,.
$$

\subsection{Why $III=\int_{\R_+^3} \text{auxiliary terms}\,dxdydt \le 0$?}
\label{III}

First of all the symmetry implies that
$$
B(u,v, \xi,\eta) = \Phi(\sqrt{u^2 +v^2}, \sqrt{\xi^2+\eta^2})\,.
$$

So far we did not use the fact that
\begin{equation}
\label{fact}
g(z) = \xi(r) +i \eta(r)\,, f(z) = e^{2i\theta} (m(r) +ik(r))\,.
\end{equation}

Let as before $a(x,y,t) = (u,v, \xi, \eta)$ with heat extension functions. Automatically, with a fixed $t$
\begin{equation}
\label{r}
\Phi (a), d\Phi (a), d^2\Phi(a)\,,\,\,\text{depend only on}\,\,\,r+\sqrt{x^2+y^2}\,.
\end{equation}

\bigskip

\noindent{\bf Remark.}
In proving that $III=0$ we are going to use this fact a lot. But $III=0$ seems to hold under some other assumptions on $f,g$.

All auxiliary terms are in
$$
\langle d^2 B(a) (z_1, z_2, \zeta_1,\zeta_2)^T, (z_2, -z_1, \zeta_2, -\zeta_1)^T\rangle -\langle d^2 B(a) (z_2, -z_1, \zeta_2,-\zeta_1)^T, (z_1, z_2, \zeta_1, \zeta_2)^T\rangle\,.
$$
This expression $= A+D_1 + D_2 +C$, where
$$
A = (B_{11}+ B_{22}) \Im z_2\bar z_1\,, C =(B_{33}+ B_{44}) \Im \zeta_2\bar \zeta_1\,.
$$
Also 
$$
D_1 = B_{13} \Im \zeta_2 \bar z_1 + B_{23} \Im \zeta_2\bar z_2 + B_{14} \Im z_1\bar \zeta_1 + B_{24} \Im z_2\bar \zeta_1\,.
$$
$$
D_2= B_{13} \Im z_2 \bar \zeta_1 + B_{23} \Im \zeta_1\bar z_1 + B_{14} \Im z_2\bar \zeta_2 + B_{24} \Im \zeta_2\bar z_1\,.
$$

\bigskip

\noindent{\bf Why $\int_{\R^2} D_1(x,y,t) \,dxdy=0$?}

In $D_1$ the smaller index of $B_{kl}$, $k\in {1,2}, l\in {3,4}$ coincides with the index of $z_i$. In $D_2$ this is not the case. This is the explanation why integrating each term of $D_1$ returns $0$. For example, (the last equality uses $\eta_{\theta}=0$)
$$
\Im \zeta_2\bar z_1 = \det\begin{bmatrix} u_x,\, \eta_x\\u_y,\, \eta_y\end{bmatrix} = \det\begin{bmatrix} u_r,\, \eta_r\\u_{\theta}/r,\, \eta_{\theta}/r\end{bmatrix}=-\eta_r u_{\theta}/r\,.
$$
But (recall $f=u+iv, g= \xi+i\eta$)
$$
B_{13} = \frac{u}{|f|}\frac{\xi}{|g|} \Phi_{12} (|f|,|g|)\,.
$$
Then the first term of $D_1$ 
$$
= \phi(r) u u_{\theta}\,,
$$
and its integral along any circle is zero. Similarly, 
\begin{equation}
\label{D1}
D_1 = (uu_{\theta} + vv_{\theta})\frac{\eta\xi_r-\xi\eta_r}{r} \frac{\Phi(|f|,|g|)}{|f||g|}\,,
\end{equation}
and so
$D_1= \phi_1(r) u u_{\theta} + \phi_2(r) vv_{\theta}$.
Hence for each fixed $t$
$$
\int_{\R^2} D_1(x,y,t) \,dxdy=0\,.
$$

\bigskip


Coming to $D_2$ we can similarly see that
\begin{equation}
\label{D2}
D_2 = (uv_{\theta}- vu_{\theta} )
 \frac{\xi\xi_r + \eta\eta_r}{r} \frac{\Phi_{12}(|f|,|g|}{|f||g|}\,.
 \end{equation}
 Recall that from \eqref{fact} it follows that $ u= m\cos 2\theta - k\sin 2\theta\,, v = m\sin 2\theta + k \sin 2\theta$, and from this
 $$
 u v_{theta} - v u_{\theta}  = 2 (m^2(r) + k^2(r)) = 2 (u^2(r) + v^2(r)) = 2 |f|^2(r) =: 2 M^2(r)\,.
 $$ 
Using similarly the notation $N(r)=|g|$ we can see  from \eqref{D2} and the previous equality that
$$
D_2= \frac2r \Phi_{12} (M(r), N(r)) N'(r) M(r)\,.
$$
Now we compute 
$$
A= (B_{11} +B_{22})\Im z_2\bar z_1 = (\Phi_{11} + \frac1{M} \Phi_1) (u_r v_{\theta}/r - v_r u_{\theta}/r)\,.
$$
Using  \eqref{fact} we get $u_r v_{\theta} - v_r u_{\theta}= 2 M(r) M'(r)$. Therefore
$$
A= \frac2r (\Phi_{11} + \frac1{M} \Phi_1)MM' = \frac2r (\Phi_{11} MM'  +\Phi_1 M')\,.
$$
Notice (again \eqref{fact}) that in 
$$
C= (B_{33} +B_{44}) \Im \zeta_2\bar \zeta_1 
$$
the expression $\Im \zeta_2\bar \zeta_1  = \xi_r \eta_{\theta}/r - \eta_r \xi_{\theta}/r =0$. So $C=0$.

Adding the expressions for $A, D_2$ we obtain after integration over $\R^2$:
\begin{align*}
\int_{\R^2} (A(x,y,t) +D_2(x,y,t)) \,dxdy = 4\pi\int_0^{\infty} (\Phi_{12}(M,N) M(r) N'(r) +\\ \Phi_{11}(M,N) M(r) M'(r)) \,dr + 4\pi\int_0^{\infty} \Phi_1(M,N) M'(r)\,dr =: a+b\,.
\end{align*}
Integrating $b$ by parts we  get $-a$ and $- \Phi_1(M(0), N(0)) M(0)$. 

\bigskip

Actually, in our particular case $III=0$. Function $f=u+iv$ on $\C$ has the form $f= e^{2i\theta} (m(r) +ik(r))$, therefore, its heat extension $f(x,y,t)$ obviously satisfies $f(x,y, 0) =0$. So $M(0)=|f(x,y,0)|=0$. As we saw
$$
III = - \Phi_1(M(0), N(0)) M(0)=0\,.
$$

\markboth{}{\sc \hfill \underline{References}\qquad}

\end{document}